\documentclass{amsart}
\usepackage{amsmath,amsthm,amssymb, amsfonts}

  \theoremstyle{definition}
  \newtheorem{theorem}{Theorem}
  \newtheorem{definition}{Definition}
  \newtheorem{example}{Example}

\title{On the permanence properties of residually exact groups}
\author{Hikaru Awazu}
\begin{document}
\maketitle

\renewcommand{\thefootnote}{\fnsymbol{footnote}}
\footnote[0]{The University of Tokyo, Graduate School of Mathematical Sciences, master's thesis.}
\renewcommand{\thefootnote}{\arabic{footnote}}

ABSTRACT.
  A discrete group $\Gamma$ is called exact if the reduced group C*-algebra ${C_{\lambda}}^{*}(\Gamma)$ is exact as C*-algebras, 
  and a discrete group $\Lambda$ is called residually exact if every nonunital element  $g \in \Lambda$ admits 
  a surjective group homomorphism from $\Lambda$ to some exact group $\Gamma$ which maps $g$ to a nonunital element of $\Gamma$.

  We prove the class of residually exact groups is closed under taking Green's graph products [1], double amalgamed products and special HNN extensions.

 \

 \

\section{Introduction}


  At first, the word ``group" tacitly denotes a countable (discrete) group in this paper.

  A discrete group is called amenable if it has a left-invariant mean. 
  This classes of groups appears in a great deal of sections of mathematics and has many characterizations with various terminologies. 
  In those, the connections with operator algebras are interesting; 
  $\Gamma$ is amenable if and only if the reduced group C*-algebra ${C_{\lambda}}^{*}(\Gamma)$ is nuclear [2]. 
  Recall that the reduced group C*-algebra ${C_{\lambda}}^{*}(\Gamma)$ is the norm closure of $\mathbb{C}[\Gamma]$ in $B(\ell^2(\Gamma))$. 

  Though the class of amenable groups and of nuclear $C^{*}$-algebras have a great deal of examples, these structures are too strong. 
  For example, the full and reduced group C*-algebras that come from an amenable group is identical.
  Then various classes weaker than those are considered, including exactness and residual amenblity.

  A C*-algebra $A$ is called exact if the functor $B \mapsto A {\otimes}_{\text{min}} B$ between the categories of C*-algebras preserves exact sequences, 
  where ${\otimes}_{\text{min}}$ is the spatial tensor product of C*-algebras.
  Exact C*-algebras can be defined another way; C*-algebras whose faithful represetations $A \hookrightarrow B(H)$ are nuclear maps [2].
  This is a natural generalizaiton of nuclerity.

  Exact groups are defined by using this exactness of C*-algebras; 
  A discrete group $\Gamma$ is called exact when its reduced group C*-algebra ${C_{\lambda}}^{*}(\Gamma)$ is exact. 

  Although the definition of exact groups comes from an operator-algebraic method, it has more group-theoretic characterizations (Yu's property A [2] and 
  Theorem 4 in this paper). 
  On the other hand, the definition of residually amenable groups is just a generalization of amenablity in the group theory, 
  imitating the construction of  residually finite groups.

\newpage


  Of course both of these have much importantance and nontrivial examples, we cannot give whole explanations in this paper. 
  For example, Baumslag-Solitar groups $B(p,q)\quad p,q \in \mathbb{N}$ are always residually amenable\footnote{At first, $BS(p,q)$ is a free-by-solvable group by [3]. 
  Recall that free groups are residually (finite $r$-) for any prime number $r$ by [9] Theorem 6.1.9 and in particular redisually solvable, 
  for finite $r$-groups are nilpotent by [9] Theorem 5.1.3. Then $BS(p,q)$ is (residually solvable)-by-amenable 
  and it becomes residually amenable by [6] Lemma 1.3.}, but not residually finite in most cases.
  Therefore the class of residually amenable groups are considerably large and hard to find non-examples.
  The class of exact groups is also large, for this includes linear groups and has various permanence properties 
  under group-operations such as taking subgroups or extensions.

  The class of residually amenable groups has been shown to be subclasses of the much larger classes of groups; sofic groups and hyperlinear groups [4][10] 
  (exact groups have not yet; see remarks of last section.). 

  However, it is remains open that all groups are sofic or hyperlinear and this problem has attracted many mathematicians.


  It is one approach to solve this problem to define new classes of groups that is ``large comparable to the sofic-class" 
  and think of the difference between them.
  This is one reason why we think of residually exact groups, a simple generalization of residual amenablity and exactness.

  Until now, there has been no study about this class and it is unknown that there are some new interesting characterizations.


  On the other hands, residual $\mathcal{C}$ had been studied by Gruenberg [5] when an arbitrary class $\mathcal{C}$ of groups is given.
  He gives a sufficient condition about $\mathcal{C}$ that the class of residual $\mathcal{C}$ is preserved by taking free products.

  Furthermore, Berlai relaxed this condition and proved the permanence properties of residually amenable groups [6] as a corolary. 
  In addition, he researched the permanence of residual $\mathcal{C}$ with respect to taking special HNN extensions, 
  double amalgamed products and graph products [7] with Ferov.

  We get use of these theorems to residual exactness and gain some permanence properties.

\section{Graph products of groups}

  Let $G = (V,E)$ be a non-oriented graph with vertexes V and edges $E$. 
  The graph $G$ is considered to be locally finite and simplical in this paper, 
  i.e. the number of edges from each vertex is finite and there is no loops and multiple edges.

  Let $\mathcal{G} = \{ {\Gamma}_v \mid v \in V \}$ be a family of groups with index set $V$.
  Then we can define the graph product $G \mathcal{G}$ of $\mathcal{G}$ with respect to $G$ as follows;

  \[
  \Gamma := *_{v\in V} {\Gamma}_v
  \]

  \[
  \Lambda := \langle[g_v,g_u] \in \Gamma| g_v \in {\Gamma}_v, g_u \in {\Gamma}_u , ({\Gamma}_v,{\Gamma}_u) \in E \rangle
  \]

  \[
  G \mathcal{G} := \Gamma / \Lambda
  \]

  where $\langle S \rangle$ denotes the normal subgroup of $G$ generated by $S \subset G$.

  Clearly it generalizes free products of groups by setting $E=\emptyset$.
  Also, if one takes a complete graph, then the graph product by this is just the direct product of groups.

\section{Residually exact groups}
\subsection{Residual class of groups in general}

  A family of groups is called class of groups if images of any groups in this family by any group isomorphisms are in the family.

  When a class $\mathcal{C}$ of groups is given, we can define residual $\mathcal{C}$ as follows.

  $\Gamma \in $residual $\mathcal{C}$ if and only if 
  for all g$\neq e \in \Gamma$, there exists $\Lambda \in \mathcal{C}$ and $\varphi : \Gamma \twoheadrightarrow \Lambda$ surjective group hom, 
  with $\varphi(g) \neq e'$.
  where $e\in \Gamma$ and $e' \in \Lambda$ are the group units.

   \

  We write Res-$\mathcal{C}$ for short.

  Note that Res-(Res-$\mathcal{C}$) = Res-$\mathcal{C}$ and for a simple group $\Gamma$, $\Gamma \in \mathcal{C} 
  \text{ if and only if } \Gamma \in \text{Res-}\mathcal{C}$.
  Moreover, Res-$\mathcal{C}$ is always closed under taking subgroups and finite direct products just by definition.

   \

  And just as an analogy of pro-finite topology, we get pro-$\mathcal{C}$ topology on a group $\Gamma$.
  More precisely, an elementary neighborhood system of the unit $e$ in pro-$\mathcal{C} $ topology on $\Gamma$ is defined to be 
  $\{\Lambda \triangleleft \Gamma \mid \Gamma/\Lambda \in \mathcal{C} \}$

\subsection{The work by Gruenberg and Berlai}

  First, we prepare some terminologies about classes of groups.

  \begin{definition}
    Let $\mathcal{C}$ and $\mathcal{D}$ be classes of groups and $\Lambda \triangleleft \Gamma$ be an arbitrary group and its normal subgroup.

    (0)  We say $\Gamma$ is a $\Lambda$-by-$\Gamma / \Lambda$ extension group in this paper.

    And we define the new class; 

    $\mathcal{C}$-by-$\mathcal{D}:=\{\Gamma \mid \Gamma \text{ is a group for which there exists } \Lambda \triangleleft \Gamma$ 
    with $\Lambda \in \mathcal{C}$ and $\Gamma / \Lambda \in \mathcal{D} \}$ 

    (1) $\mathcal{C}$ is said to satisfy Gruenberg's condition 
    if for given a normal chain of $\Gamma$; ${\Gamma}_2 \triangleleft {\Gamma}_1 \triangleleft \Gamma$ 
    with $\Gamma / {\Gamma}_1 \in \mathcal{C}$ and ${\Gamma}_1/{\Gamma}_2 \in \mathcal{C}$, 
    there exists subgroup $\Lambda \subset {\Gamma}_2$ with $\Lambda \triangleleft \Gamma$ and $\Gamma/\Lambda \in \mathcal{C}$.

     \

    (2) $\mathcal{C}$ satisfies (Gruenberg's) root property if 

    \quad (a) $\mathcal{C}$ satisfies Gruenberg's condition, 

    \quad (b) $\mathcal{C}$ is closed under taking subgroups, and 

    \quad (c) $\mathcal{C}$ is closed under taking finite direct products.

     \

    (3) A subclass $\mathcal{R} \subset \mathcal{C}$ is said to be (Berlai's) root subclass of $\mathcal{C}$ if 

    \quad (a) $\mathcal{R}$ itself ensures the root property, 

    \quad (b) $\mathcal{R}$-by-$\mathcal{C}$ $\subset$ $\mathcal{C}$ as a class of groups, and

    \quad (c) (cardinality condition) For all $\Gamma \in \mathcal{C}$, there exists $\Lambda \in \mathcal{R}$ with card($\Gamma$)$=$card($\Lambda$).
  \end{definition}

  Note that a class with Gruenberg's condition is closed under taking extensions of groups.

  And following are work by Gruenberg, Berlai and Ferov.

  \begin{theorem}
    (Gruenberg [5] Theorem4.1)

    If a class $\mathcal{C}$ of groups has the root property, then the following are equivalent;

    \quad (a) Free $\subset \text{Res-}\mathcal{C}$

    \quad (b) $\text{Res-}\mathcal{C}$ is closed under taking free products.

    where Free denotes the class of free groups and other classes of groups will be written samely.
  \end{theorem}

  \begin{theorem}
    (Berlai [6][7])

    Let $\mathcal{C}$ be a class of groups which is closed under taking finite direct products.
    And assume $\mathcal{C}$ has a root subclass $\mathcal{R}$.

    (1) Then $\text{Res-}\mathcal{C}$ is closed under taking free products.

    (2) Moreover, if $\mathcal{C}$ is closed under taking subgroups and ``Free-by-$\mathcal{C}" \subset \text{Res-}\mathcal{C}$, 

    then $\text{Res-}\mathcal{C}$ is closed under taking graph products.
  \end{theorem}

  In addition, Berlai researched special types of HNN-extensions and amalgamed free products.

  \begin{definition}
    Let $\Lambda \subset \Gamma$ a subgroup.

    (1) The HNN extension of $\Gamma$ by $id:\Lambda \rightarrow \Lambda$ is called $special$ and written as $\Gamma *_{id_{\Lambda}}$

    (2) The amalgamed product of two $\Gamma$'s by two maps incl:$\Lambda \hookrightarrow \Gamma$ is called $double$ 
    and written as $\Gamma *_{\Lambda} \Gamma$
  \end{definition}

  \begin{theorem}
    (Berlai [6])

    Let $\mathcal{C}$ be a class of groups which is closed under taking finite direct products and subgroups.
    And assume $\mathcal{C}$ has a root subclass $\mathcal{R}$.

    (1) $\Gamma *_{id_{\Lambda}} \in \text{Res-}\mathcal{C} \quad \text{ if and only if } \quad \Gamma \in \text{Res-}\mathcal{C}$ 
    and $\Lambda$ is closed in the pro-$\mathcal{C}$ topology.

    (2) Samely, $\Gamma *_{\Lambda} \Gamma \in \text{Res-}\mathcal{C} \quad \text{ if and only if } \quad \Gamma \in \text{Res-}\mathcal{C}$ 
    and $\Lambda$ is closed in the pro-$\mathcal{C}$ topology.
  \end{theorem}

\subsection{Residually exact groups}

  The definition of exact groups is shown in the introduction.
  Then we set $\mathcal{C}$ to be Exact to obtain Res-(Exact), the class of residually exact groups.

  Obviously, Res-(Finite) $\subset$ Res(Amenable) $\subset$ Res-(Exact).

  \begin{example}
    (Osajda [8])
    There is a NON-exact and residually finite group.

    Thus it is an example of NON-exact but residually exact groups.
  \end{example}

  \begin{example}
    (About Kirchberg's factorization property (F) [2].)

    The finite extensions of $SL(3,\mathbb{Z})$ (thus Finite-by-$SL(3,\mathbb{Z})$ groups) is classified completely.

    One of these, there is NON-(F) extension. But it is of course an exact group becaouse of the permanence about extensions. 
    Particularly, there is a NON-(F) but residually exact group.
  \end{example}

  \begin{example}
    (NON example of residually exactness)

    Since we have a few examples of non-exact groups, we can gain examples of non-residually finite groups from these.

    Let $\Gamma$ be a (discrete) non-exact group. 
    Because of discreteness, we can write $\Gamma = \bigcup_{F \in \mathcal{F}} F$ where $\mathcal{F} = \{ F \subset \Gamma \text{ finite set} \}$. 
    Let Alt($F$) be the alternating group on the set $F$ for $F \in \mathcal{F}$, 
    and we can regard Alt($F$) as a subgroup of Alt($F'$) where $F \subset F'$ in $\mathcal{F}$. 
    From this, we can set Alt($\Gamma$) to be the direct limit of ${(\text{Alt}(F) , \text{inclusion} )}_{F \in \mathcal{F}}$ . 

     \

    Here $\Gamma$ acts freely on Alt($\Gamma$) as follows; 

    For $\sigma \in$Alt($F$) and $g\in \Gamma$, define $g.\sigma \in$ Alt($gF$) by setting $g.\sigma (gh):= g \sigma(h)$ for $h \in F$. 
    This action is compatible for the inclusion Alt($F$) $\le$ Alt($F'$), so the action of $g$ is well-defined on whole Alt($\Gamma$). 
    Moreover, the action is free, i.e. if $g.\sigma = \sigma$ for all $\sigma \in$ Alt($\Gamma$), then $g=e$. 
    This is obvious from the definition of the action. 

     \

    Then the semidirect group Alt($\Gamma$)$\rtimes \Gamma$ becomes NON residually exact. 

    \begin{proof}
      First note that Alt($\Gamma$)$\rtimes \Gamma$ is not exact because its subgroup $\Gamma$ is not. (see $\S$4.1)

      We assume that $\varphi$:Alt($\Gamma$)$\rtimes \Gamma \twoheadrightarrow \Lambda$ is a surjective group-hom to an exact group $\Lambda$, 
      and set $e \in \Gamma$ and $e' \in \lambda$ to be units.
      Then we can show that $\varphi(\sigma,e) = e'$ even when $\sigma \neq$ identity, and thus Alt($\Gamma$)$\rtimes \Gamma$ is not residually exact.

       \

      Indeed, ker$\varphi$ $\cap$ Alt($\Gamma$) is a normal subgroup of Alt($\Gamma$) and must be trivial because Alt($\Gamma$) is simple. 
      ($\ast$1) If ker$\varphi$ $\cap$ Alt($\Gamma$) were $\{ id \}$, then we can show that $\varphi$ would be injective 
      on whole Alt($\Gamma$)$\rtimes \Gamma$. 

      Then it follows Alt($\Gamma$)$\rtimes \Gamma$ $\cong$ $\Lambda$ and contradicts with the exactness of $\Lambda$. 
      Thus ker$\varphi$ $\cap$ Alt($\Gamma$) must be Alt($\Gamma$) and $\varphi(\sigma,e) = e'$. 

       \

       To prove ($\ast$1), provided $\varphi(\sigma,g) = e'$ and $\varphi$ is injective on Alt($\Gamma$), then show $\sigma = id$ and $g = e$. 
       
        \
       (a) First we show ($\ast$2) $g.\tau = \tau$ if $\tau \sigma = \sigma \tau$ in Alt($\Gamma$). 

       $\varphi(\tau,e) = \varphi(\tau,e)\varphi(\sigma,g) = \varphi(\tau \sigma,g) 
       = \varphi(\sigma \tau, g) = \varphi(\sigma,g)\varphi(g^{-1}.\tau,e) = \varphi(g^{-1}.\tau,e). $
       
       This implies $\tau = g^{-1}.\tau$ since $\varphi$ is injective on Alt($\Gamma$). 
      	
      	\
        
       (b) Set $F_{\tau} \in \mathcal{F}$ be the set $\{h \in \Gamma \mid \tau(h) \neq h\}$ for $\tau \in$ Alt($\Gamma$). 
       \footnote{Then $F_{g.\sigma} = g F_{\sigma}$ by simple calculus and $F_{\sigma}$ is a $g$-invariant set.}
       
       Since $F_{\sigma}$ is a finite set, we can take distinct $i,j,k,l\notin F_{\sigma}$. 
       Then we can assume $i \neq gj$ and $i \neq gk$. 
       Indeed if $i=gj$, then $i = gj \neq gk,gl$, so take $(i,k,l)$ instead. 
       
      Set $\tau:= (i,j,k) \in$ Alt($\Gamma$). 
      Then ${\sigma}^{-1} \tau \sigma = (\sigma(i),\sigma(j),\sigma(k)) = (i,j,k)$, so $\sigma$ and $\tau$ commutes. 
        
      By ($\ast 2$), we obtain that $(i,j,k) = g.(i,j,k) = (gi,gj,gk)$. 
        Since we set $i \neq gj,gk$, it is shown that $i=gi$ and $g=e$. 
        
      Then $\varphi(\sigma,e) = e'$ implies $\sigma = id$ since $\varphi$ is injective on Alt($\Gamma$). 
    \end{proof}
  \end{example}

  \begin{example}
    (NON-Res-(Amenable) but Res-(Exact) group)

    Let $\Gamma$ be non-amenable but exact group. (ex. Free groups) 
    Then samely as Example.3, Alt($\Gamma)\rtimes \Gamma$ is a NON-residually-amenable but (residually) exact group. 

     \

    Indeed, Alt($\Gamma$) is exact since Exact is closed under taking direct limits of groups.
    \footnote{The reduced group C*-algebra of a direct limit of groups is just the direct limit C*algebra of their reduced group C*-algebras. 
    And exactness of C*-algebra pass to direct limits. These two fact implies the above. } 

    In addition, Alt($\Gamma)\rtimes \Gamma /$Alt($\Gamma$) $\cong$ $\Gamma$ is exact. 
    Then Alt($\Gamma)\rtimes \Gamma$ is a Exact-by-Exact group and exact itself. (See Theorem.5) 

    Simply imitating Example.3, we can prove the group is non-residually-amenable (See Theorem.10). 
  \end{example}

   \
   
   Therefore we gain the following; 

  \begin{example}

    There is NON-residually-amenable and NON-exact but residually exact group. 

    Indeed, let $\Gamma$ be a non amenable but exact group and $\Lambda$ be a non exact but residually exact group. 
    Then (Alt($\Gamma) \rtimes \Gamma) \oplus \Lambda$ is the desired one. 
  \end{example}

\section{Permanence properties of exact groups}

\subsection{Taking subgroups} 

  At first, an inclusion of discrete groups $\Lambda \hookrightarrow \Gamma$ induces the inclusion of reduced $C^*$-algebras 
  ${C_{\lambda}}^{*}(\Lambda) \hookrightarrow {C_{\lambda}}^{*}(\Gamma)$. 
  And by definition of exact $C^*$-algebras using nuclearity of a map, we easily get that C*-subalgebras of an exact C*-algebra are exact again.

  These two observations indicade that subgroups of an exact group are exact again.

\subsection{Taking direct products} \

  First we observe 

  \begin{align}
  {C_{\lambda}}^{*}({\Gamma}_1 \times {\Gamma}_2) \quad \cong \quad {C_{\lambda}}^{*}({\Gamma}_1) {\otimes}_{\text{min}} {C_{\lambda}}^{*}({\Gamma)}_2
  \end{align}

  for discrete groups ${\Gamma}_1$ and ${\Gamma}_2$.

  Immidiately we get the isomorphism of group rings; 

  \begin{align}
  \mathbb{C}[{\Gamma}_1 \times {\Gamma}_2] \quad \cong \quad \mathbb{C}[{\Gamma}_1] \odot \mathbb{C}[{\Gamma}_2]
  \end{align}

  where $\odot$ is the tensor product as vector spaces. 

  Moreover, there is the isomorphism of Hilbert spaces; 
  $\ell^2({\Gamma}_1 \times {\Gamma}_2)$ $\cong$ $\ell^2({\Gamma}_1) \hat{\otimes} \ell^2({\Gamma}_2)$
  where $\hat{\otimes}$ denotes the tensor product as Hilbert spaces.

  Then, we get $B(\ell^2({\Gamma}_1 \times {\Gamma}_2))$ $\cong$ $B(\ell^2({\Gamma}_1) \hat{\otimes} \ell^2({\Gamma}_2))$ 
  $\cong$ $B(\ell^2({\Gamma}_1)) {\otimes}_{\text{min}} B(\ell^2({\Gamma}_2))$
  ,an isometric *-isomorphism. Thus we obtain (1) taking the norm closures of (2). 

   \

  Next, the functorial definition of exact groups indicate that $A {\otimes}_{\text{min}} B$ is exact when $A$ and $B$ are exact C*-algebras. 

  Then we obtain that ${\Gamma}_1 \times {\Gamma}_2$ is exact when ${\Gamma}_1$ and ${\Gamma}_2$ is exact groups. 

\subsection{Extensions of exact groups by exact groups} \

  We need another characterization of exact groups using the concept of $amenable$ $actions$ of groups to show (Exact)-by-(Exact) $\subset$ Exact.
  Prob($\Gamma$) denotes the set of positive and norm-1 elements of $\ell^1(\Gamma)$
  for a discrete group $\Gamma$.

  \begin{definition} \

    An action of a discrete group $\Gamma$ on a topological space $X$ is called $amenable$ \
    if 
    $\text{There exists } m_i ; X \rightarrow Prob(\Gamma)$ $conti$ ,the net of maps with the index $i \in I$

    such that $\text{ for all } s \in \Gamma; \lim_{i} (\sup_{x\in X} \| s.(m_i(x)) - m_i(s.x) \|) = 0$

    where $s.m(t):= m(s^{-1}t)$ for $s,t \in \Gamma$ and $m \in Prob(\Gamma)$
  \end{definition}

  This definition extends amenablity of groups to that of actions. Indeed, actions by amenable groups must be amenable action.
  It also explains exact groups;

  \begin{theorem}
    ([2] Theorem 5.1.7) \

    For a discrete group $\Gamma$, the following are equivalent;

    (1) $\Gamma$ is exact

    (2) there exists $X$ a compact topological space and an amenable action $\Gamma \rightarrow $Homeo($X$)
  \end{theorem}

  And getting use of this characterization, we can prove the following;

  \begin{theorem}
    ([2] Prop 5.1.11)

    Let $\Lambda \triangleleft \Gamma$ be a discrete group and its normal subgroup and $X$,$Y$ be compact spaces. 
    And set $\pi$ to be the quotient map $\Gamma \twoheadrightarrow \Gamma/\Lambda$.

    Assume that there are actions $\phi: \Gamma/\Lambda \rightarrow$ Homeo($X$) and $\psi:\Gamma \rightarrow$ Homeo($Y$). 

    Moreover, if both of $\phi$ and the restricted action $\psi |_{\Lambda}$ are amenale, 

    then the map $ (\phi \circ \pi) \times \psi : \Gamma \rightarrow$ Homeo($X \times Y$) is an amenable action.
  \end{theorem}

  From this theorem, it follows that $\Gamma$ is exact if there is $\Lambda \triangleleft$ $\Lambda$ with $\Lambda$ and $\Gamma/\Lambda$ being exact. 

  Indeed, if there is an action ${\psi}':\Lambda \rightarrow$ Homeo($Y$) on just a subgroup of $\Gamma$, 
  one can get the extented action $\psi:\Gamma \rightarrow$ Homeo($Y$) and then we can get use of Theorem 5. 

  We save the proofs of these two theorems. See [2].

\section{Permanence properties of residually exact groups}

  While the class of exact groups is closed under taking extension, it is hard to confirm whether the class fulfills Gruenberg's condition or not. 
  Then we will circumvent this problem taking a nice root subclass of Berlai. 

  \begin{theorem}
  ([1])

  The class of solvable groups satisfies Gruenberg's condition and thus Gruenberg's root property.

  \end{theorem}

  We can use the class Solvable as the root subclass of Exact. 
  Indeed, because extensions of amenable groups are amenable, solvable groups must be amenable. 
  Thus Solvable $\subset$ Exact.
  Second, Solvable-by-Exact groups are exact by the corolary of Theorem 5. 
  At last, the cardinality condition (3-c) is trivial, for we now discuss countable groups. 

  In addition, Free-by-Exact groups are exact because free groups must be exact [2].

  These observations and Theorem 2 yield the following;

  \begin{theorem}
  (Main theorem) \

  The class of residually exact groups is closed under taking graph products.

  \end{theorem}

  Samely we get the results about special HNN extensions and double amalgamed products.

  \begin{theorem}
    (1) $\Gamma *_{id_{\Lambda}} \in $Res-(Exact) 
    if and only if $\Gamma \in $Res-(Exact) and $\Lambda$ is closed in the pro-Exact topology.

    (2) $\Gamma *_{\Lambda} \Gamma \in $Res-(Exact)$ \text{ if and only if } \Gamma \in $Res-(Exact) and $\Lambda$ is closed in the pro-Exact topology.
  \end{theorem}

   

\section{Acknowledgments}

  The author is very grateful to his superviser, Yasuyuki Kawahigashi for helping him on the study on master's course and proofreading this paper. 
  Also he thanks very much the expert of C*-algebra, Narutaka Ozawa for advising him on the contents of this paper.
  Moreover, he thanks his colleague Michiya Mori, Yasuhito Hashiba and Kan Kitamura for advising and encouraging.

\section{Reference}

  [1] E. Green, Graph products of groups. Ph.D. thesis, University of Leeds, 1990.

  [2] N. Brown and N.Ozawa, C*-algebras and finite-dimensional approximations. AMS, 2008.

  [3] P. H. Kropholler, Baumslag-Solitar groups and some there groups of cohomological dimension two. Comment. Math. Helvetic. {\bf 65}, 547--558, 1990.

  [4] B. Weiss, Sofic groups and dynamical systems. Sankhya Ser. {\bf A 62, 3 (2000)}, 350 –- 359. Ergodic theory and harmonic analysis, 1999.

  [5] K. W. Gruenberg, Residual properties of infinite soluble groups. Proc. London Math. Soc. {\bf (3) 7}, 1957.

  [6] F. Berlai, Residual Properties of Free Products, Communications in Algebra. {\bf 44:7}, 2959--2980, 2016.

  [7] F. Berlai and M. Ferov, Residual properties of graph products of groups. Journal of Group Theory {\bf 19(2)}, 2015.

  [8] D. Osajda, Residually finite non-exact groups. D. Geom. Funct. Anal. {\bf 28: 509}, 2018.

  [9] Derek. J. S. Robinson, A course in the theory of groups. Springer, (1982).

  [10] V. G. Pestov, Hyperlinear and sofic groups: a brief guide. Bull. Symbolic Logic {\bf 14, no. 4}, 449 -– 480, 2008.

\end{document}